\documentclass[reqno,11pt]{amsart}
\usepackage{hyperref}
\usepackage{amsmath}
\usepackage{amssymb}
\usepackage{amscd}
\usepackage{graphics}
\usepackage{latexsym}
\usepackage{stmaryrd}

\usepackage[all]{xy}

\theoremstyle{plain}
\newtheorem{theorem}{Theorem}[section]
\newtheorem{thm}[theorem]{Theorem}

\newtheorem{lem}[theorem]{Lemma}
\newtheorem{prop}[theorem]{Proposition}

\theoremstyle{definition}
\newtheorem{defn}[theorem]{Definition}

\newtheorem{ques}[theorem]{Question}
\newtheorem{rmk}[theorem]{Remark}

\newtheorem{notat}[theorem]{Notation}

\newtheorem{diag}[theorem]{Diagram}

\theoremstyle{remark}

\newcommand{\marpar}[1]{}
\newcommand{\mni}{\medskip\noindent}

\newcommand{\QQ}{\mathbb{Q}}

\newcommand{\mc}{\mathcal}

\newcommand{\pbundle}[1]{\ensuremath{\mathbb{P}({\mathcal{#1}})}}
\newcommand{\kspace}[3]{\ensuremath{\overline {\mathcal{M}}_{0,#1}(#2, #3)}}
\newcommand{\kspacenb}[3]{\ensuremath{\mathcal{M}}_{0,#1}(#2, #3)}

\newcommand{\pspace}[1]{\ensuremath{\mathbb{P}^{#1}}}

\newcommand{\Sym}{\text{Sym}}

\newsavebox{\sembox}
\newlength{\semwidth}
\newlength{\boxwidth}

\newsavebox{\semrbox}
\newlength{\semrwidth}
\newlength{\boxrwidth}

\title[Families of Rational Curves]
{Complete families of linearly non-degenerate rational curves}

\author[DeLand]{M\lowercase{atthew} D\lowercase{e}L\lowercase{and}}
\address{Department of Mathematics \\
  Columbia University\\ New York, NY 10025}
\email{deland@math.columbia.edu} 

\date{\today}
\begin{document}


\begin{abstract}
We prove that a complete family of linearly non-degenerate rational curves of degree $e > 2$ in $\pspace{n}$ has at most $n-1$ moduli.  For $e = 2$ we prove that such a family has at most $n$ moduli.  It is unknown whether or not this is the best possible result.  The general method involves exhibiting a map from the base of a family $X$ to the Grassmaninian of $e$-planes in $\pspace{n}$ and analyzing the resulting map on cohomology.
\end{abstract}


\maketitle

\section{Introduction and Main Theorem}

\mni
Let $Y$ be a smooth, projective variety over $\mathbb{C}$.  The Kontsevich moduli space $\kspace{0}{Y}{\beta}$ parametrizes isomorphism classes of pairs $(C,f)$ where $C$ is a proper, connected, at-worst-nodal, arithmetic genus $0$ curve, and $f$ is a stable morphism $f:C\rightarrow Y$ such that $f_*[C] = \beta \in H_2(Y,\mathbb{Z})$.  This is a Deligne-Mumford stack whose coarse moduli space, $\overline{M}_{0,0}(Y,\beta)$ is projective.  See, for example, [FP].  

\mni
For the remainder of this paper, we will restrict to the case of degree $e$ curves in $Y = \pspace{n}$.

\mni
Let $\mathcal{U} \subset \kspacenb{0}{\pspace{n}}{e}$ be the open substack parametrizing maps $f: \pspace{1} \rightarrow \pspace{n}$ which are isomorphisms onto their image, such that the span of each image is a $\pspace{e}$.  Note that no point in $\mc{U}$ admits automorphisms, and that $\mathcal{U}$ is isomorphic to an open subscheme in the appropriate Hilbert and Chow schemes. In particular, $\mc{U}$ itself is a scheme.

\begin{defn} \label{main_def}
 Suppose $X$ and $\mc{C}$ are proper varieties and $\pi: \mc{C} \rightarrow X$ is a proper surjective morphism.  We will consider diagrams of the form:
$$
\xymatrix{
\mathcal{C}  \ar[r]^{f} \ar[d]^{\pi} & \mathbb{P}^n \\
   X }
$$

\mni
In the case where each fiber of $\pi$ is a $\pspace{1}$ and $f$, restricted to each fiber, corresponds to a point in $\mc{U}$, we will call the diagram a \textit{complete family of linearly non-degenerate degree $e$ curves}.  Such a family induces a map $\alpha: X \rightarrow \mc{U}$.  If the map is generically finite, that is, if $\dim X = \dim \alpha(X)$ we will call the diagram a \textit{family of maximal moduli}.  We will refer to $X$ as the base of the family.  Note that $\mc{C}$ is the pullback of the universal curve over $\mc{U}$, and so we will refer to the map $f$ above as $ev$.  The notation $(\mc{C}, X, ev, \pi, n, e)$ will denote a complete family of linearly non-degenerate degree $e$ curves in $\pspace{n}$, as above.
\end{defn}

\mni
One can ask for the largest number of moduli of such a family, that is, the dimension of the base $X$ of a family of maximal moduli.  The bend and break lemma [DEB] gives a strict upper bound on the dimension of complete subvarieties of $X \in \kspacenb{0}{\pspace{n}}{e}$, namely $2n - 2$.  M. Chang and Z. Ran prove that if $\Lambda$ is an effectively parametrized family of curves in $\pspace{n}$ then $\dim \Lambda \leq n - 1$ [CR].  Note that a complete family of linearly non-degenerate curves is not effectively parametrized in the sense of [CR] because each fiber is a $\pspace{1}$.  Also, a simple corollary to the theorem of Coskun-Harris-Starr shows that the number of moduli of a linearly non-degenerate family of degree $e$ curves in $\pspace{e}$ is in fact $0$ [CHS].  The main result of this paper is:

\begin{thm} \label{thm-main}
If $X$ is the base of a family of linearly non-degenerate degree $e \geq 3$ curves in $\pspace{n}$ with maximal moduli, then $\dim X \leq n - 1$.  If $X$ is the base of such a family of non-degenerate degree $2$ curves in $\pspace{n}$, then $\dim X \leq n$.
\end{thm} 

\begin{rmk}
It is unclear that this is the best possible result.  There are certainly examples of $r$ dimensional families in $\pspace{r+e}$.  One way to construct such families is to take the Segre embedding:
$$
\xymatrix{
\pspace{1} \times \pspace{r} \ar[r]^<<<<{(e,1)} & \mathbb{P}^N}
$$

\noindent where $N = (e + 1) \cdot (r + 1) - 1$.  Project from a point $p \in \pspace{N}$ not in any $\pspace{e}$ spanned by the image of $\pspace{1}\times \{q\}$ for every point $q \in \pspace{r}$.  This gives an $r$ dimensional family of non-degenerate degree $e$ curves in $\pspace{N-1}$.  Continue projecting in this fashion.  We can always find a point $p$ to project from as long as $N > r + e$.  So we arrive at an $r$ dimensional family of degree $e$ curves in $\pspace{r + e}$. 

\end{rmk}

\subsection{Discussion}

\begin{ques} 
Can an $r$ dimensional family of degree $e$ non-degenerate rational curves be constructed in $\pspace{m}$ for $m < r+ e$?
\end{ques}
\begin{ques}
Our bound is obviously not optimal when $e < n$ and by the theorem of [CHS] mentioned above, nor when $n = e$ either.  Other small dimensional examples remain unknown to me.  Is it possible to have a $2$ dimensional family of smooth conics in $\pspace{3}$ or a $2$ dimensional family of smooth cubics in $\pspace{4}$?
\end{ques}
\begin{ques}
If the variety swept out by these curves is required to be contained in  a smooth hypersurface, does the bound improve?  In fact, this question was the original motivation for this work.
\end{ques}

\subsection{Outline of Proof}

\mni
We give a brief outline of the proof:

\mni
Let $e > 2$ and fix $X$ to be the base of a complete family of linearly non-degenerate degree $e$ curves in $\pspace{n}$ with maximal moduli. Assume that $\dim X \geq n$.  By results from section 2, we will reduce the situation to the case where the universal curve $\mathcal{C}$ over $X$ is the projectivization of a rank $2$ vector bundle $\mathcal{E}$ on $X$.  The situation will then be further reduced to the case where we have the following maps:

\begin{diag} \label{diag}
$$
\xymatrix{
& \mathcal{C} =  \pbundle{E} \ar[r]^{ev} \ar[d]^{\pi} & \mathbb{P}^n &\\
  \pbundle{E} \ar[r] \ar@/_/[ddrr]_\gamma & X \ar@/^/[dr]^\phi  & &\\
   & & Gr(e+1,n+1) & Fl(1, \ldots, n) \ar[l] \ar[dl] \\
   & & Fl(1, \ldots, e + 1) \ar[u]&}
$$
\end{diag}

\mni
where $\phi$ and $\gamma$ are generically finite morphisms and the two maps from $Fl(1, \ldots, n)$ are the natural projection maps.

\mni
In section 3, we will construct an ample line bundle $\mc{L}$ on $Fl(1, \ldots, e + 1)$ and give a cohomology argument to show that $c_1(\mc{L})^{n + 1}$ pulls back to $0$ under $\gamma$. This will allow us to conclude.  In the case $e = 2$, a different computation is needed, but similar ideas apply.

\begin{notat}
Fix the ambient $\pspace{n}$.  We will denote by $Fl(a_1, \ldots, a_k)$ with $a_1 < a_2 < \ldots < a_k$ the flag variety parameterizing vector subspaces $A_k \subset A_{k-1}  \subset \ldots \subset A_{1} \subset \mathbb{C}^{n+1}$ such that $\text{codim} (A_i, \mathbb{C}^{n+1}) = a_i$.  In the special case $Fl(a)$ we will write $Gr(a, n+1)$, the set of $a$ dimensional quotients of $\mathbb{C}^{n+1}$.  Also $\pbundle{E}$ will refer to the set of hyperplanes in the fibers of $\mathcal{E}$.  This is the convention used in [EGA II], confusingly dual to the one used in [HAR] and [FUL], all of which are references for this paper.
\end{notat}

\mni
I happily thank my advisor, Aise Johan de Jong, for many helpful discussions, suggestions, and for his untiring patience. 

\section{Reductions}
\mni
We will first prove some general lemmas which will soon be applied to the case of a complete family of linear non-degenerate degree $e$ curves.




\begin{prop} \label{projectivize}
Suppose that $\pi: \mc{C} \rightarrow X$ is a proper, surjective morphism of complete varieties where each fiber of $\pi$ is abstractly isomorphic to $\pspace{1}$. Then there exists a generically finite map $f:X' \rightarrow X$ such that in the fiber square:
$$
\xymatrix{
\mathcal{C}' \ar[d]^{\pi'} \ar[r]^{f'} & \mathcal{C}  \ar[d]^{\pi} \\
   X' \ar[r]^f & X }
$$
\noindent
$\pi'$ realizes $\mathcal{C}'$ as the projectivization of a rank $2$ vector bundle $\mc{E}$ on $X'$.  That is, $\mathcal{C}' = \pbundle{E}$.
\end{prop}

\begin{proof}
Let $\mu$ denote the generic point of $X$ and $k(X)$ denote the function field.  Each fiber of $\pi$ is just a $\pspace{1}$, so $\pi$ is a smooth morphism [HAR, III.10].  From this, it follows that the relative dualizing sheaf $\mathcal{L} = \omega_{ \mathcal{C}/X } |_{\pi^{-1}(\mu)}$ is isomorphic to the canonical sheaf on the generic fiber $\mathcal{C}_\mu$ [ibid].  We then use $\mathcal{L}^{-1}$ to embed $\mathcal{C}_\mu$ as a degree $2$ rational curve in $\pspace{2}_{k(X)}$.  The image may not have a $k(X)$ point, but after a base change to a degree two field extension $K \supset k(X)$ we can arrange that $(\mathcal{C}_\mu)_K \subset \pspace{2}_K$ will have a $K$ point.

\mni
Let $X'$ be the normalization of $X$ in $K$ [EGA II, 6.3] and pull back $\mathcal{C}$ to $X'$.  We explain the following diagram:

$$
\xymatrix{
\mathcal{C'} \ar[r] \ar[d] & X' \ar@{-->}@/_/[l] \ar[d]_<<<<f & Spec(K) \ar[d] \ar[l] \ar[dll]|<<<<<<<<<<<<<<<\hole \\
\mathcal{C} \ar[r] & X & Spec( k(X) ) \ar[l]}
$$

\mni
The left square is defined as the fiber product of $X'$ and $\mc{C}$ over $X$.  The map $f: X' \rightarrow X$ is normalization map, which is finite.  By the properties of normalization, the function field of $X'$ is $K$, that is, $X'$ has a $K$ point [ibid].  And since $\mathcal{C}$ has a $K$ point, there is a rational map from $X'$ to $\mathcal{C}$ inducing a rational map from $X'$ to $\mathcal{C}'$, the dotted arrow above.  This rational map can be resoloved by blowing up $X'$ along some ideal sheaf [HAR].  That is, we have

$$
\xymatrix{
\mathcal{C}'' \ar[d] \ar[r] & X'' \ar@/_1pc/[l]_\sigma \ar[d]^g \ar[dl]^h \\
\mathcal{C}' \ar[d] \ar[r] & X' \ar[d]^f \\
\mathcal{C}  \ar[r] & X}
$$

\mni
The top square is defined so that $\mathcal{C}''$ is the fiber product.  Then the maps $h,g$ determine a section $\sigma: X'' \rightarrow \mathcal{C}''$.  Of course, since $f$ is finite, and $g$ is generically finite, the map $f \circ g$ is also generically finite.  Note that each geometric fiber of $\mathcal{C}'$ over $X'$ is still a $\pspace{1}$, and that the same holds for geometric fibers of $\pi'': \mathcal{C}'' \rightarrow X''$.  The existence of the section will allow us to conclude that $\mc{C}'' \cong \pbundle{E}$ by a standard argument:

\mni
Let $D = \sigma(X'')$ be a divisor on $\mathcal{C''}$, so $D.F = 1$ where
$F$ is the class of a fibre, and hence $\mc{L}(D)$ is relatively very ample.  By Grauert's Theorem (HAR III.12), $\mathcal{E} = \pi''_*\mathcal{L}(D)$ is locally free of rank $2$ on $X$.  The natural map from $\pi''^* \mathcal{E}$ to $\mathcal{L}(D)$ is surjective:  We can check this fiberwise by Nakayama's Lemma.  Each fiber $\mathcal{C}''_x$ is just $\pspace{1}$, $\mathcal{L}(D)$ is generated by global sections on each fiber, and $\mathcal{E} \otimes k(x) \rightarrow H^0(\mathcal{L}(D)_x)$ is surjective by Grauert's Theorem again.  This surjection $\pi''^* \mathcal{E} \rightarrow \mathcal{L}(D)$ determines a morphism $\beta:\mathcal{C}'' \rightarrow \pbundle{E}$ (HAR) such that $\beta^*\mathcal{O}(1) \simeq \mathcal{L}(D)$. Now $\beta$ is an isomorphism on each fiber, so an isomorphism.

\mni
Thus constructing the section was enough to show that we have the following picture:

$$
\xymatrix{
\mathcal{C}'' \ar[d]^{\pi''} \ar[r]^{f'} & \mathcal{C}  \ar[d]^\pi \\
   X' \ar@/^/[u]^\sigma \ar[r]^f & X }
$$

\mni
where $f$ is generically finite, and $\mathcal{C}'' \simeq \pbundle{E}$ for a rank $2$ vector bundle $\mathcal{E}$ on $X'$.  This completes the proof.
\end{proof}



\mni
In the case where a projective bundle over $X$ admits a map to $\pspace{n}$, we would like to say something about the pullback of $\mc{O}_{\pspace{n}}(1)$:

\begin{prop} \label{flatten_bundle}
Suppose that $\mc{E}$ is a rank $2$ vector bundle on a variety $X$ and let $\pi: \pbundle{E} \rightarrow X$ be the natural map.  Suppose in addition that $\pbundle{E}$ admits a map to $\pspace{n}$ which is degree $e$ on each fiber. Then there exists a finite map $f: X' \rightarrow X$ such that in the fiber product diagram:

$$
\xymatrix{
\pbundle{E_{\text{X}'}} \ar[d]^{\pi'} \ar[r]^{f'} & \pbundle{E}  \ar[r]^{ev} \ar[d]^{\pi} & \mathbb{P}^n \\
   X' \ar[r]^f & X }
$$

\mni
we have that $\pi'_* ev'^* \mathcal{O}(1) = \Sym^e(\mathcal{E_{\text{X}'}})$ where $ev' = ev \circ f'$.
\end{prop}

\begin{proof}
First we remark that $ev^*\mathcal{O}(1)$ is a line bundle that is degree $e$ on each fiber of $\pi$.  Thus $ev^*\mathcal{O}(1) = \mathcal{O}(e) \otimes \pi^*(N)$ for some line bundle $N$ on X. This follows by the description of the Picard group of a projective bundle [HAR].  Then $\pi_* ev^*\mathcal{O}(1) = \Sym^e(\mathcal{E}) \otimes \mc{N}$.  If there is a line bundle $\mathcal{L}$ on $X$ such that $\mathcal{L}^e \simeq \mc{N}$ then it is an easy exercise to show that $\Sym^e(\mathcal{E}) \otimes N \simeq \Sym^e(\mathcal{E} \otimes \mathcal{L})$ and it is well known [HAR] that $\pbundle{E} \simeq \pbundle{E \otimes L}$.  The following Lemma then allows us to conclude.
\end{proof}

\begin{lem} \label{finite_extension}
Let $\mc{N}$ be a line bundle on a variety $X$.  There is a finite map $\tau: X' \rightarrow X$ and a line bundle $\mc{L}$ on $X'$ such that $\mc{L}^{\otimes e} \simeq  \mc{N}_\tau$, where $\mc{N}_\tau$ denotes the pullback of $\mc{N}$ to $X'$.  
\end{lem}
\begin{proof}
Let $\pi: \mc{N} \rightarrow X$ be the structure map.  Choose an open affine covering of $X$, ${U_i = Spec(A_i)}$ where $\mc{N}$ is trivialized.  That is, for each $i$ we have:
$$
\xymatrix{
\pi^{-1}(U_i) \ar[r]^{\phi_i} \ar[d]^\pi & U_i \times \mathbb{A}^1 \ar[dl]^\pi \\
U_i &}
$$

\mni
where $\phi_i$ is an isomorphism and $\pi$ is actually the restiction of $\pi$ to $\pi^{-1}(U_i)$.  On overlaps $U_{ij} = U_i \cap U_j$, we have isomorphisms $\phi_i^{-1} \circ \phi_j: U_{ij} \times \mathbb{A}^1 \rightarrow U_{ij} \times \mathbb{A}^1$ given by a global section $r_{ij}$ of $\mc{O}_{U_{ij}}$, that is, a rational function on $X$.  These elements $r_{ij}$ are subject to the usual cocycle conditions: 
\begin{align*}
r_{ii} &= 1 \text{ in } A_i \\
r_{ij} r_{ji} &= 1 \text { in } A_{ij} \\
r_{ki} r_{ij} r_{jk} &= 1 \text { in } A_{ijk}
\end{align*}

\mni
where $A_{ij}$ (respecively $A_{ijk}$) is $\mc{O}(U_{ij})$ (respectively $\mc{O}(U_{ijk})$).  Let $k$ denote the function field of $X$.  We will adjoin an $e^{\text{th}}$ root of each $r_{ij}$, call it $s_{ij}$ to $k$ and argue that this can be done consistently.  That is, the $s_{ij}$ can be chosen to satisfy the analogous cocycle conditions stated above.  The idea is that having chosen, for example $s_{12}$ and $s_{13}$, the section $s_{23}$ is determined automatically by the third cocycle condition.  Then choosing $s_{14}$ now forces the choice of $s_{24}$ and $s_{34}$.  Then one checks that the choice of $s_{24}$ satisfies all cocycle relations determining it at this point.  Continuing in this fashion, the choices of $s_{1k}$ determine every other $e^{\text{th}}$ root.  Adjoin each of these elements $s_{ij}$ to $k$ and call $L$ the resulting finite field extension.  Let $B_i$ (respectively $B_{ij}$) be the integral closure of $A_i$ (respectively $A_{ij}$) in $L$.  Note that $s_{ij}$ is contained in $B_{ij}$.  Finally, let $X'$ be the scheme where $Spec(B_i)$ are glued together using the $Spec(B_{ij})$.  Let the line  bundle $\mc{L}$ be determined by the $s_{ij}$.  The map $\tau:X' \rightarrow X$ is finite because $B_i$ is module finite over $A_i$.  By construction, $\mc{L}^e = \mc{N}_\tau$ on $X'$.
\end{proof}



\section{Proof}

\mni
Before looking at the general case, we first prove a stronger result than the main theorem would imply when $n=e$:
 
\begin{prop} \label{prop-base-case}
If $n = e$, and $(\mc{C}, X, ev, \pi, n, n)$ is a family of maximal moduli as in Definition~\ref{main_def}, then $\dim X = 0$.  That is, there is no complete curve contained in $\mc{U} \subset \kspacenb{0}{\pspace{n}}{n}$.
\end{prop}

\begin{proof}
Suppose that there is such a family with $\dim X > 0$.  We apply a result of Coskun, Harris, and Starr where they have computed the effective cone of $\kspace{0}{\pspace{d}}{d}$.  By Theorem 1.5 in [CHS] the effective cone of this space is generated by non-negative linear combinations of the divisor class $\mathcal{D}_{deg}$ and classes $\Delta_{ij}$ supported on the boundary.  Recall from their paper that $\mathcal{D}_{deg}$ denotes the locus of maps where the linear span of the image is not the entirety of $\pspace{d}$.  As before, let $\alpha(X)$ be the image of $X$ in $\mc{U}$.  By the properties of $X$, all divisors coming from the boundary restrict to $0$ on $\alpha(X)$, as does $\mathcal{D}_{deg}$ because $X$ is the base of a linearly non-degenerate family.  If the entire effective cone restricts to $0$ on $\alpha(X)$, then certainly the ample cone does as well.  This is a contradiction, because if $Y$ is a complete variety, and if $Y'$ is a complete, irreducible, subvariety of $Y$, then an ample divisor on $Y$ restricts to have positive degree on $Y'$: see, for example [HAR, app A].  Thus $\dim X = \dim \alpha(X) = 0$.
\end{proof}

\mni
We are now ready to prove the main theorem.
\begin{proof}[Proof of Theorem ~\ref{thm-main}]

\mni
Fix $(\mc{C}, X, ev, \pi, n, e)$ to be a family of maximal moduli as in Defintion~\ref{main_def} with $2 < e < n$.  By way of contradiction, assume that $\dim X \geq n$.  By taking an irreducible proper subvariety of $X$, and restricting the family, we may assume that $\dim X = n$.

\mni
For any point $x \in X$, denote by $\phi(x)$ the linear $e$-plane spanned by the image of the map corresponding to $x$.  That is, $\phi(x) = \text{Span}(ev(\pi^{-1}(x))$.  The map $\phi: X \rightarrow Gr(e+1, n+1)$ is well-defined because each curve corresponding to a point in $X$ is linearly non-degenerate.  The morphism is quasifinite by Proposition \ref{prop-base-case} but it is proper, so finite.

\mni
Applying Proposition~\ref{projectivize} and then Proposition~\ref{flatten_bundle} we may assume that there is a generically finite map $f: X' \rightarrow X$ such that we have fiber product diagram:

$$
\xymatrix{
\pbundle{E} \ar[d]^{\pi'} \ar[r]^{f'} & \mc{C}  \ar[r]^{ev} \ar[d]^{\pi} & \mathbb{P}^n \\
   X' \ar[r]^f & X }
$$

\mni
where $\mc{E}$ is a rank two vector bundle on $X'$ and $\pi'_*(f' \circ ev)^* \mc{O}(1) = \Sym^e(\mc{E})$. The collection $(\pbundle{E}, X', f' \circ ev, \pi', n, e)$ is still a family of linearly non-degenerate degree $e$ curves with maximal moduli, and $\dim X' = n$.  The composed map $f \circ \phi$ is a generically finite map from $X'$ to the Grassmanian.  To simplify notation, we rename this new family $(\pbundle{E}, X, ev, \pi, n, e)$ and trust that no confusion will arise.

\mni
We construct the universal section.  Let $Y = \pbundle{E}$ and consider the fiber product diagram:

$$
\xymatrix{
\pbundle{E_{\text{Y}}} \ar[d]^{\pi'} \ar[r] & \pbundle{E}  \ar[d]^{\pi} \\
   Y \ar[r] & X }
$$

\mni
We have a natural section $\sigma: Y \rightarrow \pbundle{E_{\text{Y}}}$ given by the diagonal map.  This section corresponds to a surjection $\mc{E}_Y \rightarrow \mathcal{L}$ where $\mathcal{L} = \sigma^*\mathcal{O}_{\pbundle{E_\text{Y}}}(1)$.  Let $\mathcal{L}_1 = \mathcal{L}$ and let $\mathcal{L}_2$ be the line bundle such that:

$$
0 \rightarrow \mathcal{L}_2 \rightarrow \mathcal{E_\text{Y}} \rightarrow \mathcal{L}_1 \rightarrow 0
$$

\mni
This sequence induces a filtration on $\Sym^e(\mathcal{E})$:

$$
\Sym^e( \mathcal{E_\text{Y}} ) = F^0 \supset F^1 \supset \dots F^e \supset F^{e+1} = 0
$$

\mni
such that $F^p/F^{p+1} \simeq \mathcal{L}_2^p \otimes \mathcal{L}_1^{e-p}$ (HAR II.5).  Note that $Y$ corresponds to curves parameterized by $X$ and a point on that curve.  We have a natural map from $Y \rightarrow Gr(e+1,n+1)$ by composition, and the data of the $F^p$s induce a map from $\gamma: Y \rightarrow Fl(1, \ldots, e+1)$.  Informally, the information of ``the point" on the curve induces a linear filtration of the $\pspace{e}$ spanned by the curve.  The linear spaces in between the point and the entire $\pspace{e}$ are the osculating $k$-planes, $k = 1, \ldots, e$. We can see this by working locally where the map is defined by $t \rightarrow (1, t, t^2, \ldots, t^e, 0, \ldots, 0)$. All the maps in diagram \ref{diag} have been constructed.

\mni
On $Fl(1, \ldots, e+1)$ we have the natural sequence of universal quotient bundles:

$$
\mathcal{O}^{n+1} \rightarrow \mathcal{Q}_{e+1} \rightarrow \dots \rightarrow \mathcal{Q}_1 \rightarrow 0
$$

\mni
Recall the previously constructed map: $\gamma: \pbundle{E} \rightarrow Fl(1,\ldots,e + 1)$.  The proof hinges on the fact that we can construct an ample line bundle on the flag manifold whose first chern class pulled back to $\pbundle{E}$ raised to the $(n+1)^{st}$ power is zero.

\mni
For $p = 0, \ldots, e$ let $x_p = c_1( \ker \mathcal{Q}_{p+1} \rightarrow \mathcal{Q}_p)$.  By construction of $\gamma$ we have $\gamma^* x_p = c_1(F_p / F_{p+1}) = pc_1(\mathcal{L}_2)+(e-p)c_1(\mathcal{L}_1)$.

\mni
Consider the projection map $pr : Fl(1, \ldots, n) \rightarrow Fl(1, \ldots, e+1)$ and the injective map it induces on cohomology (always with rational coefficients):

$$
pr^*: H^*( Fl(1, \ldots, e+1) ) \rightarrow H^*( Fl(1,\ldots,n) )
$$

\mni 
It is well known that $H^*( Fl(1,\ldots,n) ) = \QQ[x_0,\ldots,x_n]/\mathcal{I}$ where $\mathcal{I}$ is the ideal of symmetric polynomials in the $x_i$s [FUL].  By a slight abuse of notation, call $pr^* (x_i) = x_i$.

\mni
In the cohomology ring of full flags, we claim that $x_p^{n+1} = 0$ for each $p$.  To see this, note that in this ring, the following identity holds:

$$
T^{n+1} = (T - x_1) \cdot (T - x_2) \cdot \dots \cdot (T- x_n)
$$

\mni
since on the right hand side each coefficient of $T^k$ with $k < n+1$ is a symmetric polynomial.  Taking $T = x_p$ proves the identity.  Then since $pr^*$ is injective, we must also have that $x_p^{n+1} = 0$ in the cohomology ring of partial flags, so $(pc_1(\mathcal{L}_2)+(e-p)c_1(\mathcal{L}_1))^{n+1} = 0$ for each $p = 0, \ldots, e$.

\mni
To simplify notation, in what follows we write $z = c_1(\mathcal{L}_1)$ and $y = c_1(\mathcal{L}_2)$. For relevant facts about the cohomology ring of the flag variety, see Appendix~\ref{sec_flag}.  For any $D = \lambda_0 x_0 + \ldots + \lambda_{e} x_{e}$  we have:

\begin{align*}
\gamma^*(D) &= \gamma^*(\lambda_0 \cdot x_0 + \ldots + \lambda_{e} \cdot x_{e})\\
&= \sum_{p=0}^{e} \lambda_p \cdot (py + (e-p)z) \\
&= (\lambda_1 + 2\lambda_2 + 3\lambda_3 + \ldots  + e\lambda_e)y + (e\lambda_0 + (e-1)\lambda_1 + \ldots + \lambda_{e-1})z 
\end{align*}

\mni
Let $A$ be the coefficient of $y$ and $B$ the coefficient of $z$.  If we can choose $\lambda_0, \ldots, \lambda_e$ so that $\gamma^*(D) = Ay + Bz$ is a $\QQ$ multiple of one of the $(py + (e-p)z)$ then for some rational number $m$ we have:

\begin{align*}
\gamma^*( D^{n+1}) &= (m (py + (e-p)z))^{n+1} \\
&= 0
\end{align*}

\mni
Suppose that we can also arrange that $D$ is the chern class of an ample line bundle on $Fl(1, \ldots, e + 1)$.  Then we have [FUL] that $D^{\dim Y} \cdot \gamma(Y) > 0$ because $\gamma$ is generically finite.  Since $\dim Y = n + 1$ we can rewrite this as $ (D|_{\gamma(Y)})^{n+1} > 0$.  Applying Lemma~\ref{proj_fmla}, we see that $\gamma^*(D)^{n+1} > 0$ which contradicts the above calculation.  Hence we can conclude that $\text{dim} \pbundle{E} < n + 1$ and so $\text{dim} X < n$.

\mni
It remains to show that $D$ can be chosen with these properties.  See Appendix~\ref{sec_flag} for a description of the ample cone of the flag variety.  To arrange this choice of $D$, set 

$$
\lambda_0 = \frac{1}{e}, \lambda_1 = \frac{1}{e-1}, \ldots, \lambda_i = \frac{1}{e-i}, \ldots, \lambda_{e-1} = 1.
$$

\mni
Then obviously we have that $B = e$.  We will prove that $\lambda_e$ can be chosen to satisfy:

\begin{align*}
\lambda_e > \lambda_{e-1} = 1 \text{ and } \frac{A}{B} & = e - 1 \\
\intertext{This is equivalent to:} 
 e\lambda_e & = e(e-1) - \Sigma_{i=1}^{e-1} \frac{i}{e-i} \\
 \lambda_e & = (e - 1) - \Sigma_{i=1}^{e-1} \frac{i}{e(e-i)} >^? 1 \\
\intertext{Using partial fractions and simplifying, we get} \\
 \lambda_e & = e - \Sigma_{i = 0}^{e-1} \frac{1}{e - i} \\
\end{align*}

\mni
It is then easy to show this is strictly larger than $1$ as long as $e \geq 3$.  Therefore $D$ can be chosen with the required positivity property and the proof is complete when $e \geq 3$.  A simple calculation shows this method cannot work when $e=2$.  To show a slightly weaker result in that case, we need another method.
\end{proof}

\mni
We include the statement of the projection formula used in the proof above:

\begin{lem} \label{proj_fmla}
[DEB] Let $\pi :  V \rightarrow W$ be a surjective morphism between proper varieties.  Let $D_1, \ldots, D_r$ be Cartier divisors on $W$ with $r \geq \dim(V)$.  Then the projection formula holds, i.e.:
$$
\pi^*D_1 \cdots \pi^*D_r = \deg(\pi)(D_1 \cdots D_r)
$$
\end{lem}

\section{The Proof for Conics}

\mni
In this section we prove a bound for families of smooth conics one dimension weaker than for a family of higher degree curves.  Note that for conics (and cubics), being linearly non-degenerate is equivalent to having smooth images.

\begin{thm}
If $(\mc{C}, X, ev, \pi, 2, n)$ is a family of linearly non-degenerate conics in $\pspace{n}$ with maximal moduli, then $\dim X \leq n$.
\end{thm}
\begin{proof}
Exactly as in the case $e > 2$, we apply Proposition~\ref{projectivize} and then Proposition~\ref{flatten_bundle} to reduce to the case where the family has the form:

$$
\xymatrix{
\mc{C} = \pbundle{E} \ar[d]^{\pi} \ar[r]^{ev} & \mathbb{P}^n \\
   X &}
$$

\mni
where $\mc{E}$ is a rank two vector bundle on $X$ and $\pi_* ev^* \mc{O}(1) = \Sym^2(\mc{E})$.  As in the higher degree case, we have a generically finite map $\phi: X \rightarrow Gr(3, n + 1)$.  On the Grassmanian $Gr(3, n+1)$, we have the tautological exact sequence:

$$
0 \rightarrow \mc{S} \rightarrow \mc{O} \rightarrow \mc{Q} \rightarrow 0
$$

\mni
where $\mc{Q}$ is the tautolocial rank $3$ quotient bundle.  Applying Lemma~\ref{finite_extension} and pulling back the family one more time, we may further assume that $\phi^*(\mc{Q}) = \Sym^2(\mc{E})$.  

\mni
Now we proceed with a Chern class computation.  First, we compute the Chern polynomial:

\begin{align*}
c_t(\Sym^2(\mc{E})) &= 1 + 3c_1(\mc{E})t + (2c_1(\mc{E})^2 + 4c_2(\mc{E}))t^2 + 4c_1(\mc{E})c_2(\mc{E})t^3
\end{align*} 

\mni
If we let $A = 3c_1(\mc{E})$, $B = 2c_1(\mc{E})^2 + 4c_2(\mc{E})$, and $C = 4c_1(\mc{E})c_2(\mc{E})$, an easy computation shows
$$
9AB - 27C - 2A^3 = 0
$$
Write $\widetilde{A} = c_1(\mc{Q})$, $\widetilde{B} = c_2(\mc{Q})$, and $\widetilde{C} = c_3(\mc{Q})$.  These classes pull back under $\phi$ in the following way:

$$
A = c_1(\Sym^2(\mc{E})) = c_1(\phi^*(\mc{Q})) = \phi^*(c_1(\mc{Q})) = \phi^*(\widetilde{A})
$$

\mni
Here, we have used the properties of $\phi$ and the functoriality of Chern classes.  Similarly $B = \phi^*(\widetilde{B})$ and $C = \phi^*(\widetilde{C})$.  By the functoriality of Chern classes and the above relationships, we have

$$
\phi^*(9\widetilde{A}\widetilde{B} - 27\widetilde{C} -2\widetilde{A}^3) = 0
$$

\mni
Let $\xi = 9\widetilde{A}\widetilde{B} - 27\widetilde{C} -2\widetilde{A}^3$.  It becomes convenient to rewrite $\xi$ in terms of the chern roots of $\mc{Q}$.  If $\alpha_1, \alpha_2, \alpha_3$ are the Chern roots of $\mc{Q}$, then we calculate:

\begin{align*}
\widetilde{A} &=  \alpha_1 +\alpha_2 + \alpha_3 \\
\widetilde{B} &= \alpha_1\alpha2 + \alpha_1\alpha_3 + \alpha_2\alpha_3 \\
\widetilde{C} &= \alpha_1\alpha_2\alpha_3 \\
\xi &= (\alpha_1 + \alpha_2 - 2\alpha_3)(\alpha_2 + \alpha_3 - 2\alpha_1)(\alpha_1 + \alpha_3 - 2\alpha_2)
\end{align*}

\mni
Now let $f = \phi_*[X] \in H^*(Gr(3, n+1), \mathbb{Q} )$ where $[X]$ is the fundamental class of $X$.  The projection formula then gives $\xi \cdot f = 0$.  

\mni
Since $c_1(\mc{Q})$ is positive, $c_1(\phi^*\mc{Q})$ is positive by Lemma~\ref{proj_fmla}, and we get the desired bound on $\dim X$ by showing that $c_1(\phi^*\mc{Q})^{n+1} = 0$.  Since we have already shown that $\phi^*(\xi) = 0$, it would suffice to show that $c_1(\mc{Q})^{n+1}$ is divisible by $\xi$ in $H^*(Gr(3,n+1))$.  Instead, we show that this relationship holds in the cohomology ring of full flags, and argue that this is enough to conclude.

\mni
\textit{Claim}: $\xi$ divides $(\alpha_1 + \alpha_2 + \alpha_3)^{n+1}$ in $H^*(Fl, \mathbb{Q} )$, where $Fl$ denotes the space of full flags.

\mni
Consider the following fiber square:

$$
\xymatrix{
\widetilde{X} \ar[d]^{p'} \ar[r]^{\phi'} & Fl \ar[d]^p \\
X \ar[r]^{\phi} & Gr(3, n+1)
}
$$ 

\mni
We have presentations for the cohomology rings:

$$ 
H^*(Gr, \mathbb{Q} ) = \mathbb{Q}[\alpha_1, \alpha_2, \alpha_3] / I
$$
$$
H^*(Fl, \mathbb{Q} ) = \mathbb{Q}[\alpha_1, \ldots, \alpha_{n+1}] / (Symm)
$$

\mni
where $Symm$ is the ideal generated by the elementary symmetric functions, and the injective map $p^*$ satisfies $p^*(\alpha_i) = \alpha_i$ for $i = 1,2,3$.  In $H^*(Fl, \mathbb{Q} )$ we have 

$$
T^{n+1} = (T - \alpha_1) \cdots (T - \alpha_{n+1})
$$

\mni
as before.  Evaluate the two sides of the equation at $T = \frac{\alpha_1 + \alpha_2 + \alpha_3}{3}$ to find:

\begin{align*}
(\alpha_1 + \alpha_2 + \alpha_3)^{n+1} &= (\frac{\alpha_2 + \alpha_3 - 2\alpha_1}{3})(\frac{\alpha_1 + \alpha_3 - 2\alpha_2}{3})(\frac{\alpha_1 + \alpha_2 - 2\alpha_3}{3})g'(\alpha) \\
&= \xi \cdot g(\alpha)
\end{align*}

\mni
for some polynomials $g'$ and $g$ which proves the claim.  To finish the proof, remark that the fibers of $p$ are projective varieties, that is, effective cycles, and so the same is true of $p'$.  By [FUL], we have

$$
(p')^* \phi^* (c_1(\mc{Q}))^{n+1} = (\phi')^*p^*(c_1(\mc{Q}))^{n+1}
$$

\mni
The left hand side of the equation gives an effective cycle on $\widetilde{X}$, in particular, a non-zero cohomology class.  On the right side, however, we get:

\begin{align*}
(\phi')^*p^*(c_1(\mc{Q}))^{n+1} &= (\phi')^*(\alpha_1 + \alpha_2 + \alpha_3)^{n+1} \\
&= (\phi')^*(\xi \cdot g(\alpha)) \\
&= (\phi')^*( p^* \xi \cdot g(\alpha)) \\
&= (\phi')^*p^*\xi \cdot (\phi')^* g(\alpha) \\
&= (p')^*\phi^*\xi \cdot (\phi')^* g(\alpha) \\
&= 0 \cdot (\phi')^* g(\alpha) \\
&= 0
\end{align*}

\mni
This gives a contradiction, so we conclude that $\dim(X) \leq n$.
\end{proof}

\section{Appendix - Divisors on the Flag Variety} \label{sec_flag}
\mni
We include some notes on the ample cone of the flag variety $F = Fl(1,\ldots, e+1)$.  Let $w_i$ be the $\pspace{1}$ constructed by letting the $i^\text{th}$ flag vary while leaving the others constant.  These $e+1$ lines freely generate the homology group $H_2(F)$ and the effective cone of curves.  The $e+1$ chern classes $x_p = c_1(ker(\mc{Q}_{p+1} \rightarrow \mc{Q}_p))$ generate $H^2(F)$ and we check that the intersection matrix $\langle x_i, w_j \rangle$ is given by :

$$
\begin{pmatrix}
1 & 0 & \ldots & 0 &0 \\
-1 & 1 & \ldots &0 & 0 \\
\vdots & \vdots& \ddots &\vdots & \vdots \\
0 & 0 & \ldots & 1 & 0 \\
0 & 0 & \ldots & -1 & 1\\
\end{pmatrix}
$$

\mni
with $1$'s on the diagonal and $-1$'s on the lower diagonal.  The ample cone of $F$ is given by combinations of the $x_i$'s which evaluate positively.  That is, by $\QQ$ divisors $\lambda_0 x_0 + \ldots + \lambda_{e} x_{e}$ where $\lambda_0 > 0$, $\lambda_1  > \lambda_0$, $\ldots$ , $\lambda_{e} > \lambda_{e-1}$. 

\mni
In fact, it is well known that for varieties of the type $F = G/B$, the Picard Group of $F$ is isomorphic to the character group of $F$, often denoted $X(T)$ where $T$ is a maximal torus.  Any character can be written as a linear combination of the fundamental weights $\lambda = \sum a_i t_i$ and a character is called dominant if all $a_i \geq 0$, regular if all $a_i$ are non-zero.  The ample divisors correspond exactly to the dominant and regular characters.  [LG].  In our case, the full flag variety corresponds to $G/B$ for $G = SL(n + 1)$.  The simple roots correspond to $s_i = \alpha_i - \alpha_{i+1}$ for $0 \leq i \leq n$.  Suppose $L = \lambda_1 x_0 + \ldots + \lambda_{n} x_{n}$ where the $x_i$ are as above.  Then $L$ corresponds to the weight $\lambda_0 s_0 + \ldots + \lambda_{n} s_{n}$ which is dominant if and only if $L$ is ample, if and only if $\lambda_1 > 0$, $\lambda_2 > \lambda_1$, $\ldots$, $\lambda_{n+1} > \lambda_{n}$.  The case of the partial flag variety then follows immediately from this one.


\begin{thebibliography}{9}
\setlength{\baselineskip}{12pt}

\bibitem[CHS]{CHS} 
I.~Coskun, J.~Harris., and J.~Starr.
\newblock{The Effective cone of the Kontsevich Moduli Space}.
\newblock{To Appear.}
\newblock{See http://www.math.sunysb.edu/~jstarr/papers/index.html}.

\bibitem[CR]{changran} 
M.~Chang and Z.~Ran.
\newblock{Closed Families of Smooth Space Curves}.
\newblock{\em Duke Mathematical Journal} {\bf 52}(1985), no 3, 707-713.

\bibitem[DEB]{Deb} 
O.~Debarre.
\newblock{\em Higher Dimensional Algebraic Geometry}.
\newblock Springer-Verlag, New York, 2001.
\newblock Universitext.

\bibitem[EGA]{EGA}
A.~Grothendieck and J~Dieudonn\'{e}.
\newblock{\em El\'{e}ments de G\'{e}om\'{e}trie Alg\'{e}brique. II. \'{E}tude globale \'{e}l\'{e}mentaire de quelques classes de morphismes}
\newblock{Publ. Math. IHES} {\bf 8}(1961), 5-222.
  

\bibitem[FP]{fultonpan:kontsevich}
W.~Fulton and R.~Pandharipande.
\newblock {Notes on stable maps and quantum cohomology}.
\newblock In {\em Algebraic geometry---Santa Cruz 1995}, volume 62 Part 2 of
  {\em Proc. Sympos. Pure Math.}, pages 45--96. Amer. Math. Soc., 1997.


\bibitem[FUL]{fulton:intersection}
W.~Fulton.
\newblock {\em Intersection theory}, volume~2 of {\em Ergebnisse der Mathematik
  und ihrer Grenzgebiete. 3. Folge. A Series of Modern Surveys in Mathematics}.
\newblock Springer-Verlag, Berlin, second edition, 1998.


\bibitem[HAR]{hartshorne:book}
R.~Hartshorne.
\newblock {\em Algebraic geometry}.
\newblock Springer-Verlag, New York, 1977.
\newblock Graduate Texts in Mathematics, No. 52.

\bibitem[LG]{LG} 
V.~Lakshamibai and N~Gonciulea.
\newblock {\em Flag Varieties}
\newblock Hermann, Editeurs Des Sciences et Des Arts, 2001.
\newblock Travaux en Cours, No. 63.
 
\end{thebibliography}
\end{document}